\documentclass[a4paper,11pt]{article}

  \usepackage[english]{babel}
\usepackage{amsfonts,amsmath,amsthm,graphicx,a4wide}    
\usepackage{graphicx} 
\usepackage{amsmath} 
\usepackage{amssymb}  
\usepackage[colorlinks=true]{hyperref}
\usepackage[utf8]{inputenc}
\usepackage{mathrsfs}
\usepackage{dsfont} 
\usepackage{xcolor}

\newtheorem{theorem}{Theorem}[section]

\newtheorem{lemma}[theorem]{Lemma}
\newtheorem{proposition}[theorem]{Proposition}

\newtheorem{definition}[theorem]{Definition}
\newtheorem{remark}[theorem]{Remark}

\title{\LARGE \bf
On exponential stabilization of two-qubit systems
}

\author{W. Liang, N. H. Amini, P. Mason
\thanks{All authors are with Laboratoire des Signaux et Syst\`emes, CNRS - CentraleSup\'elec - Univ. Paris-Sud, Universit\'e Paris-Saclay, 3, rue Joliot Curie, 91192, Gif-sur-Yvette, France. {\tt\small [first name].[family name]@l2s.centralesupelec.fr.}This work was supported by the Agence Nationale de la Recherche project QUACO ANR-17-CE40-0007}}

\date{}
\begin{document}

\maketitle

\begin{abstract}
In this paper, we consider a two-qubit system undergoing
continuous-time measurements. In presence of multiple channels, we provide sufficient
conditions on the continuous feedback control law ensuring almost sure exponential convergence to a predetermined Bell state. This is obtained by applying
stochastic tools, Lyapunov methods and geometric control tools.  With one channel, we establish asymptotic convergence towards a predetermined Bell state. In both cases, we provide explicit expressions of feedback control laws satisfying the above-mentioned conditions. Finally, we demonstrate the
effectiveness of our methodology through numerical simulations.
\end{abstract}
\section{Introduction}
In view of the rapid development of quantum information science~\cite{nielsen2010quantum}, the generation of quantum entangled states~\cite{bengtsson2017geometry} has become essential  in a variety of applications such as quantum teleportation, quantum cryptography and quantum computation. The simplest entangled states are the Bell states, which are pure states corresponding to maximal quantum entanglement of two spin-$\frac12$ systems {(i.e. a two-qubit system)}. 

Two-qubit systems undergoing continuous-time measurements  represent a particular example of open quantum systems whose evolution is described by stochastic master equations.  The problem of controlling open quantum systems by feedback  has received attention in the control community mainly starting from the years 2000s~\cite{van2005feedback, armen2002adaptive, mirrahimi2007stabilizing, tsumura2008global,ahn2002continuous,yamamoto2007feedback,mabuchi2005principles}. This field is a branch of the stochastic control theory whose first main ideas have been developed by Belavkin in~\cite{belavkin1983theory}.

Concerning stabilization of two-qubit systems, some interesting results have been derived in~\cite{mirrahimi2007stabilizing} and~\cite{yamamoto2007feedback}.   In~\cite{mirrahimi2007stabilizing}, the authors design a switching quantum feedback controller that asymptotically stabilizes the  system towards  two specific Bell states. Then, in~\cite{yamamoto2007feedback}, the methods in~\cite{van2005feedback} are adopted in order to construct a continuous feedback controller stabilizing  the target Bell state starting from almost any initial pure state  when the measurement is perfect.

 In this paper we generalize the methods in~\cite{liang2018exponential,liang2019exponential} in order to apply them to the feedback stabilization problem for two-qubit systems with multiple quantum channels. In this case, unlike the above-mentioned papers, the associated measurement operators may be degenerate rendering the stabilization problem more complicated.  Here, we  derive some general conditions on the feedback laws and the control Hamiltonians enforcing the exponential convergence towards the target Bell state. In addition, when only one quantum channel is available, we design a continuous feedback law which asymptotically stabilizes the system to an arbitrary Bell state.
This control law is adapted from the  switching feedback law  proposed in~\cite{mirrahimi2007stabilizing}, although our stability analysis  is radically different.  

This paper is organized as follows. In Section~\ref{SEC:SYS}, we introduce the stochastic model describing two-qubit systems with multiple quantum channels in presence of imperfect measurements. In Section~\ref{SEC:PRE}, we introduce the notions of stochastic stability needed throughout the paper. In Section~\ref{SEC:2QC}, we show the exponential convergence of the system with two quantum channels and with zero control input towards the set of Bell states. We consider  a class of appropriate  control Hamiltonians and further propose necessary conditions on the continuous feedback law ensuring the exponential stabilization. In Section~\ref{SEC:1QC}, explicit feedback laws are exhibited which asymptotically stabilize the system with only one quantum channel to the target state. Simulation results are provided in Section~\ref{SEC:SIM}.

\textit{Notations:}
The imaginary unit is denoted by $i$. We take $\mathds{1}$ as the identity matrix. We denote the conjugate transpose of a matrix $A$ by $A^*.$ $A_{i,j}$ represents the element of the matrix $A$ at $i$-th row and $j$-th column. The function $\mathrm{Tr}(A)$ corresponds to the trace of a square matrix $A.$ The commutator of two square matrices $A$ and $B$ is denoted by $[A,B]:=AB-BA.$ $A\otimes B$ represents the Kronecker product of $A$ and $B$. For $x\in\mathbb{C}$, $\mathbf{Re}\{x\}$ is the real part of $x$ and $\mathbf{Im}\{x\}$ is the imaginary part of $x$. 
\section{System description}
\label{SEC:SYS}
The dynamics of a two-qubit system undergoing continuous-time measurements with $n$ channels can be described by a matrix-valued stochastic differential equation of the form~\cite{van2005feedback,hudson1984quantum,belavkin1989nondemolition,bouten2007introduction}:
\begin{equation}
d\rho_t=F_0(\rho_t)dt+\sum^n_{k=1}F_{k}(\rho_t)dt+\sum^n_{k=1}\sqrt{\eta_{k}}G_k(\rho_t)dW_{k}(t),
\label{SME}
\end{equation}
where the quantum state is described by the density operator $\rho_t$, which belongs to the compact space $\mathcal{S}:=\{\rho\in\mathbb{C}^{4 \times 4}|\,\rho=\rho^*,\mathrm{Tr}(\rho)=1,\rho\geq0\}$. 
Here $W_t=(W_k(t))_{1\leq k\leq n}$ is a $n$-dimensional standard Wiener process with the natural filtration $\mathcal{F}_t,$  and $W_k$ are independent, i.e., for $i,j=1,\dots,n,$ one has $\left\langle W_i(t),W_j(t)\right\rangle=\delta_{i,j}t.$ The filtered probability space associated with the above evolution is $(\Omega,\mathcal{F},(\mathcal{F}_t),\mathbb{P}).$ The measurement efficiency for the  $k$-th channel  is given by $\eta_k\in(0,1]$. The functions $F_0$, $F_k$ and $G_k$ are given by the following expressions
\begin{equation*}
\begin{split}
F_0(\rho)&:=-i[H_0,\rho]-i\textstyle\sum_{j=1}^m u_j(\rho)[H_j,\rho],\\
F_k(\rho)&:=L_k\rho {L_k}-L^2_k\rho/2-\rho L^2_k/2,\\
G_k(\rho)&:=L_k\rho+\rho L_k-2\mathrm{Tr}(L_k\rho)\rho.
\end{split}
\end{equation*}
{The function $u$ appearing in $F_0$}
denotes the feedback law taking values in $\mathbb{R}^m,$ while $H_0=H^*_0\in\mathbb{C}^{4 \times 4}$ is the free Hamiltonian, $H_j=H^*_j\in\mathbb{C}^{4 \times 4}$ are the control Hamiltonians and  $L_k=L^*_k\in\mathbb{C}^{4 \times 4}$ are the measurement operators. If the feedback $u$ is in $\mathcal{C}^1(\mathcal{S},\mathbb{R}^m)$, the existence and uniqueness of the solution of~\eqref{SME} as well as the strong Markov property of the solution are ensured by the results established in~\cite{mirrahimi2007stabilizing}. 

In the following sections, we consider the feedback stabilization of the above system towards one of the four Bell states $\boldsymbol\Psi_{\pm}=\Psi_{\pm}\Psi^*_{\pm}$ and $\boldsymbol\Phi_{\pm}=\Phi_{\pm}\Phi^*_{\pm}$ where 
{
\begin{align*}
\Psi_{\pm}=
\frac{1}{\sqrt{2}}\left(
\left[\begin{matrix}1 \\ 0 \end{matrix}\right]\otimes \left[\begin{matrix}1 \\ 0 \end{matrix}\right] 
\pm \left[\begin{matrix}0 \\ 1 \end{matrix}\right] \otimes \left[\begin{matrix}0 \\ 1 \end{matrix}\right]\right) = \frac{1}{\sqrt{2}}\left[\begin{smallmatrix}1 \\ 0 \\ 0\\ \pm 1\end{smallmatrix}\right],\\
\Phi_{\pm}=
\frac{1}{\sqrt{2}}\left(
\left[\begin{matrix}1 \\ 0 \end{matrix}\right]\otimes \left[\begin{matrix}0 \\ 1 \end{matrix}\right] 
\pm \left[\begin{matrix}0 \\ 1 \end{matrix}\right] \otimes \left[\begin{matrix}1 \\ 0 \end{matrix}\right]\right) = \frac{1}{\sqrt{2}}\left[\begin{smallmatrix}0 \\ 1 \\ \pm 1\\0\end{smallmatrix}\right].
\end{align*}
}
These states are assumed to be common eigenstates of the measurement operators $L_k$. 
\section{Preliminary stochastic tools}
\label{SEC:PRE}
In this section, we introduce some basic definitions which are fundamental for the rest of the paper.
\paragraph{Stochastic stability}
We introduce some notions of stochastic stability needed throughout the paper by adapting classical notions (see e.g.~\cite{mao2007stochastic, khasminskii2011stochastic}) to our setting. In order to provide them, we first present the definition of Bures distance~\cite{bengtsson2017geometry}.
\begin{definition}
The Bures distance between two quantum states $\rho_a$ and $\rho_b$ in $\mathcal{S}$ is defined as
\begin{equation*}
\setlength{\abovedisplayskip}{3pt}
\setlength{\belowdisplayskip}{3pt}
d_B(\rho_a,\rho_b) := \sqrt{2-2\mathrm{Tr}\left( \sqrt{\sqrt{\rho_b}\rho_a\sqrt{\rho_b}} \right)}.
\end{equation*}
Also, the Bures distance between a quantum state $\rho_a$ and a set $E \subseteq \mathcal{S}$ is defined as 
\begin{equation*}
\setlength{\abovedisplayskip}{3pt}
\setlength{\belowdisplayskip}{3pt}
d_B(\rho_a, E): = \min_{\rho \in E} d_B(\rho_a,\rho).
\end{equation*}
\end{definition}
Given $E \subseteq \mathcal{S}$ and $r>0$, we define the neighborhood $B_r(E)$ of $E$ as
\begin{equation*}
\setlength{\abovedisplayskip}{3pt}
\setlength{\belowdisplayskip}{3pt}
B_r(E) := \{\rho \in \mathcal{S}|\, d_B(\rho,E) < r\}.
\end{equation*}
\begin{definition}
Let $\bar E$ be an a.s. invariant set of system~\eqref{SME}, then $\bar E$ is said to be
\begin{enumerate}
\item[1.] 
\emph{locally stable in probability}, if for every $\varepsilon \in (0,1)$ and for every $r >0$, there exists  $\delta = \delta(\varepsilon,r)$ such that,
\begin{equation*}
\setlength{\abovedisplayskip}{3pt}
\setlength{\belowdisplayskip}{3pt}
\mathbb{P} \left( \rho_t \in B_r (\bar E) \text{ for } t \geq 0 \right) \geq 1-\varepsilon,
\end{equation*}
whenever $\rho_0 \in B_{\delta} (\bar E)$.

\item[2.]
\emph{almost surely asymptotically stable}, if it is locally stable in probability and,
\begin{equation*}
\setlength{\abovedisplayskip}{3pt}
\setlength{\belowdisplayskip}{3pt}
\mathbb{P} \left( \lim_{t\rightarrow\infty}d_B(\rho_t,\bar E)=0 \right) = 1,
\end{equation*}
whenever $\rho_0 \in \mathcal{S}$.

\item[3.]
\emph{exponentially stable in mean}, if for some positive constants $\alpha$ and $\beta$,
\begin{equation*}
\setlength{\abovedisplayskip}{3pt}
\setlength{\belowdisplayskip}{3pt}
\mathbb{E}(d_B(\rho_t,\bar E)) \leq \alpha \,d_B(\rho_0,\bar E)e^{-\beta t},
\end{equation*} 
whenever $\rho_0 \in \mathcal{S}$. The smallest value $-\beta$ for which the above inequality is satisfied is called the \emph{average Lyapunov exponent}. 

\item[4.]
\emph{almost surely exponentially stable}, if
\begin{equation*}
\setlength{\abovedisplayskip}{3pt}
\setlength{\belowdisplayskip}{3pt}
\limsup_{t \rightarrow \infty} \frac{1}{t} \log d_B(\rho_t,\bar E) < 0, \quad a.s.
\end{equation*}
whenever $\rho_0 \in \mathcal{S}$. The left-hand side of the above inequality is called the \emph{sample Lyapunov exponent} of the solution. 
\end{enumerate}
\end{definition}
 Note that any equilibrium $\bar{\rho}$ of~\eqref{SME}, that is any quantum state satisfying $\sum^n_{k=0}F_k(\bar{\rho})=0$ and $G_k(\bar{\rho})=0$ for $k\in\{1,\dots,n\}$, is a special case of invariant set. 
\paragraph{Infinitesimal generator and It\^o's formula}
Given a stochastic differential equation $dq_t=f(q_t)dt+g(q_t)dW_t$, where $q_t$ takes values in $Q\subset \mathbb{R}^p,$ the infinitesimal generator  is the operator $\mathscr{L}$ acting on twice continuously differentiable functions $V: Q \times \mathbb{R}_+ \rightarrow \mathbb{R}$ in the following way
\begin{equation*}
\setlength{\abovedisplayskip}{3pt}
\setlength{\belowdisplayskip}{3pt}
\begin{split}
\mathscr{L}V(q,t):=\frac{\partial V(q,t)}{\partial t}&+\sum_{i=1}^p\frac{\partial V(q,t)}{\partial q_i}f_i(q)+\frac12 \mathrm{Tr}\left(g^*(q)\frac{\partial^2 V(q,t)}{\partial q_i\partial q_j}g(q)\right).
\end{split}
\end{equation*}
It\^o's formula describes the variation of the function $V$ along solutions of the stochastic differential equation and is given as follows
$dV(q,t) = \mathscr{L}V(q,t)dt+\frac{\partial V(q,t)}{\partial q_i}g(q)dW_t.$

\medskip

With these basic tools, we are ready to analyse the stabilization problem of Equation~\eqref{SME}.  From now on, the operator $\mathscr{L}$ is associated with the equation~\eqref{SME} and we choose appropriate measurement operators so that Bell states are equilibria of~\eqref{SME}. In particular, we set $\bar{E}:=\{\boldsymbol\Psi_{\pm},\boldsymbol\Phi_{\pm}\}$. 
\section{System with two quantum channels}
\label{SEC:2QC}
In this section, we study the dynamics of the system~\eqref{SME} with two quantum channels, i.e., with two measurement operators $L_1$ and $L_2$ and, consequently, two diffusion terms. 
{Consider the Pauli matrices
\begin{equation*}
\setlength{\abovedisplayskip}{3pt}
\setlength{\belowdisplayskip}{3pt}
\sigma_x = \left(\begin{matrix}0 &1\\ 1 &0 \end{matrix}\right),\ \sigma_y=\left(\begin{matrix}0 &-i\\ i &0 \end{matrix}\right),\ \sigma_z=\left(\begin{matrix}1 &0\\ 0 &-1 \end{matrix}\right).
\end{equation*}}
We set $L_1 = \sqrt{M_1}L_z$ with $M_1>0$ and $L_2 = \sqrt{M_2}L_x$ with $M_2>0$, where $L_z:=\sigma_z\otimes\sigma_z$ and $L_x:=\sigma_x\otimes\sigma_x$. Here $M_1, M_2>0$ are the strengths of the interaction between the light and the atoms. We also take $H_0=\omega L_z$ with $\omega\geq 0$ and use only one control Hamiltonian $H_1$.
{Note that the four Bell states coincide with the common eigenstates of the chosen operators $L_1$ and $L_2$.}
\subsection{Quantum State Reduction}
Before attacking the exponential stabilization problem of the system~\eqref{SME}, we analyze the asymptotic behavior of the system~\eqref{SME} with $u(\rho)\equiv 0$. This highlights the importance of the diffusion terms to design the feedback laws as already observed in~\cite{liang2018exponential,liang2019exponential,abe2008analysis}. 

We define $\Lambda_{1}(\rho):=\rho_{1,1}+\rho_{4,4}$, $\Lambda_{2}(\rho):=\rho_{2,2}+\rho_{3,3}$, $\mathcal{S}_{\Psi} := \{\rho\in\mathcal{S}|\,\Lambda_{2}(\rho)=0\},$  $\mathcal{S}_{\Phi} := \{\rho\in\mathcal{S}|\,\Lambda_{1}(\rho)=0\}$ and $\Delta(\rho):=\Lambda_{1}(\rho)\mathbf{Re}\{\rho_{2,3}\}-\Lambda_{2}(\rho)\mathbf{Re}\{\rho_{1,4}\}$. Then we state the following two lemmas, inspired by analogous result in~\cite{mao2007stochastic,khasminskii2011stochastic}, identifying invariant subsets of the system. 
\begin{lemma}
Assume $u\equiv 0$. Then $\mathcal{S}_{\Psi}$ and $\mathcal{S}_{\Phi}$ are a.s. invariant for Equation~\eqref{SME}. If the initial state satisfies $\rho_0\notin \mathcal{S}_{\Psi}$ or $\rho_0\notin \mathcal{S}_{\Phi}$, then $\mathbb{P}( \rho_t\notin \mathcal{S}_{\Psi}, \forall\, t\geq 0 )=1$ or $\mathbb{P}( \rho_t\notin \mathcal{S}_{\Phi}, \forall\, t\geq 0 )=1$ respectively.
\label{Never reach without feedback 2_1} 
\end{lemma}
\proof
{
For $u\equiv 0,$ the dynamics of $\Lambda_{2}(\rho)$ 
is given by 
\begin{equation*}
\setlength{\abovedisplayskip}{3pt}
\setlength{\belowdisplayskip}{3pt}
\begin{split}
d\Lambda_{2}(\rho_t)&=-4\sqrt{\eta_1 M_1}\Lambda_{1}(\rho_t)\Lambda_{2}(\rho_t)dW_1(t)+4\sqrt{\eta_2M_2}\Delta(\rho_t)dW_2(t).
\end{split}
\end{equation*}
Since $\rho\geq 0$ one has 
\begin{equation*}
\setlength{\abovedisplayskip}{3pt}
\setlength{\belowdisplayskip}{3pt}
|\mathbf{Re}\{\rho_{2,3}\}|\leq \Lambda_{2}(\rho)/2,\quad |\mathbf{Re}\{\rho_{1,4}\}|\leq \Lambda_{1}(\rho)/2,
\end{equation*}
so that $\Delta(\rho)=0$ for all $\rho\in\mathcal{S}_{\Psi}\cup\mathcal{S}_{\Phi}$, which yields the first part of the lemma for $\mathcal{S}_{\Psi}$.

Let us now prove the second part of the lemma. Given $\varepsilon>0$, consider any $\mathcal{C}^2$ function on $\mathcal{S}$ such that 
\begin{equation*}
\setlength{\abovedisplayskip}{3pt}
\setlength{\belowdisplayskip}{3pt}
\begin{split}
V_{\Psi}(\rho)&=1/\Lambda_{2}(\rho),\quad \mbox{if }\Lambda_{2}(\rho)>\varepsilon.
\end{split}
\end{equation*}
A simple computation shows that $\mathscr{L}V_{\Psi}(\rho)\leq RV_{\Psi}(\rho)$ if $\Lambda_{2}(\rho)>\varepsilon$ 
 for some positive constant $R$.
To conclude the proof in the case $\rho_0\notin \mathcal{S}_{\Psi}$, one just applies the same arguments as in~\cite[Lemma 4.1]{liang2019exponential}. 
Roughly speaking, setting $f(\rho,t)=e^{-R t} V_{\Psi}(\rho)$ one has $\mathscr{L}f\leq 0$ whenever $\Lambda_{2}>\varepsilon$. From this fact one proves that the probability of reaching $\mathcal{S}_{\Psi}$ in a finite fixed time $T$ is proportional to $\varepsilon$ and, being the latter arbitrary, it must be $0$. This conclude the proof in the case $\rho_0\notin \mathcal{S}_{\Psi}$. The same arguments lead to the result in the case $\rho_0\notin \mathcal{S}_{\Phi}$.}
\hfill$\square$

We denote $\Gamma(\rho):=2\mathbf{Re}\{\rho_{2,3}\}+2\mathbf{Re}\{\rho_{1,4}\}$, $\mathcal{R}_{+}:=\{\rho\in\mathcal{S}|\,\Gamma(\rho)=1\}$ and $\mathcal{R}_{-}:=\{\rho\in\mathcal{S}|\,\Gamma(\rho)=-1\}.$
\begin{lemma}
Assume $u\equiv 0$.  If $\rho_0\in\mathcal{S}_{\Psi}\setminus\{\boldsymbol\Psi_{\pm}\}$ or $\rho_0\in\mathcal{S}_{\Phi}\setminus\{\boldsymbol\Phi_{\pm}\}$, then $\mathbb{P}( \rho_t\neq\boldsymbol\Psi_{\pm}, \forall\, t\geq 0 )=1$ or $\mathbb{P}( \rho_t\neq\boldsymbol\Phi_{\pm}, \forall\, t\geq 0 )=1$ respectively.
\label{Never reach without feedback 2_2} 
\end{lemma}
\proof
As $u\equiv 0,$ then $\boldsymbol\Psi_{\pm}$ and $\boldsymbol\Phi_{\pm}$ are four equilibria of Equation~\eqref{SME}. Moreover, we have the following relations
\begin{equation*}
\setlength{\abovedisplayskip}{3pt}
\setlength{\belowdisplayskip}{3pt}
\mathcal{S}_{\Psi}\cap\mathcal{R}_{\pm}=\boldsymbol\Psi_{\pm},\quad \mathcal{S}_{\Phi}\cap\mathcal{R}_{\pm}=\boldsymbol\Phi_{\pm}. 
\end{equation*}
The dynamic of $\Gamma(\rho)$ is given by 
\begin{equation*}
\setlength{\abovedisplayskip}{3pt}
\setlength{\belowdisplayskip}{3pt}
\begin{split}
d\Gamma(\rho_t)&=2\sqrt{\eta_2M_2}(1-\Gamma^2(\rho_t))dW_2(t)-8\sqrt{\eta_1M_1}\Delta(\rho_t)dW_1(t).
\end{split}
\end{equation*}
Due to Lemma~\ref{Never reach without feedback 2_1}, for all $\rho\in\mathcal{S}_{\Psi}\cup\mathcal{S}_{\Phi}$, we have
\begin{equation*}
\setlength{\abovedisplayskip}{3pt}
\setlength{\belowdisplayskip}{3pt}
d\Gamma(\rho_t)=2\sqrt{\eta_2M_2}(1-\Gamma^2(\rho_t))dW_2(t).
\end{equation*}
Given $\varepsilon>0$, consider any $\mathcal{C}^2$ functions on $\mathcal{S}$ such that 
\begin{equation*}
\begin{split}
V_+(\rho)&=1/(1+\Gamma(\rho)),\quad \mbox{if }1+\Gamma(\rho)>\varepsilon;\\
V_-(\rho)&=1/(1-\Gamma(\rho)),\quad \mbox{if }1-\Gamma(\rho)>\varepsilon.
\end{split}
\end{equation*}
Then we have $\mathscr{L}V_+(\rho)\leq R_+V_+(\rho)$ if $1+\Gamma(\rho)>\varepsilon$ and $\mathscr{L}V_-(\rho)\leq R_-V_-(\rho)$ if $1-\Gamma(\rho)>\varepsilon$ for some positive constants $R_+$ and $R_-.$ To conclude the proof, one just applies the same arguments as in {the previous lemma and}~\cite[Lemma 4.1]{liang2019exponential}. \hfill$\square$

For the system~\eqref{SME}, the asymptotic convergence towards $\bar{E}=\{\boldsymbol\Psi_{\pm},\boldsymbol\Phi_{\pm}\}$ has been proved in~\cite{vu2012real}. 
{We now show the exponential convergence towards $\bar E$ in mean and almost surely.}
\begin{theorem}[Quantum state reduction]
For system~\eqref{SME}, with $u\equiv0$ and $\rho_0 \in \mathcal{S},$ the set $\bar{E}$ is exponentially stable in mean and a.s. with average and sample Lyapunov exponent less or equal than $-\min\{\eta_1 M_1,\eta_2M_2\}$. Moreover, the probability of convergence to $\bar{\boldsymbol\rho}\in\bar{E}$ is $\mathrm{Tr}(\rho_0 \bar{\boldsymbol\rho})$.
\label{QSR_2}
\end{theorem}
\proof
Consider the function
\begin{equation}
V(\rho)=\sqrt{V_z(\rho)+V_x(\rho)}
\label{LYA_QSR_2}
\end{equation}
as a candidate Lyapunov function, where $V_z(\rho):=\Lambda_{1}(\rho)\Lambda_{2}(\rho)$ and $V_x(\rho):=1-\Gamma^2(\rho)$. Note that $4V_z(\rho)=\mathrm{Tr}(L^2_z\rho)-\mathrm{Tr}^2(L_z\rho)$ and $V_x(\rho)=\mathrm{Tr}(L^2_x\rho)-\mathrm{Tr}^2(L_x\rho)$ can be considered as the ``variance" functions of $L_z$ and $L_x$ respectively. Moreover, $V(\rho)=0$ if and only if $\rho\in\bar{E}$. Due to Lemmas~\ref{Never reach without feedback 2_1} and \ref{Never reach without feedback 2_2}, {the complementary set of $\bar E$ is invariant. Since $V$ is twice continuously differentiable in this set, we can apply It\^o's formula. We have} $\mathscr{L}V(\rho) \leq -\bar{C}V(\rho)$ with $\bar{C}:=\min\{\eta_1 M_1,\eta_2 M_2\}$.
For all $\rho_0 \in \mathcal{S}$, we have 
\begin{equation*}
\setlength{\abovedisplayskip}{3pt}
\setlength{\belowdisplayskip}{3pt}
\mathbb{E}(V(\rho_t)) = V(\rho_0)-\bar{C}\int^t_0 \mathbb{E}(V(\rho_s))ds.
\end{equation*}
In virtue of Gr\"onwall inequality, we have
$
\mathbb{E}(V(\rho_t))\leq V(\rho_0) e^{-\bar{C}t}.
$
By a straightforward calculation, we can show that the candidate Lyapunov function is bounded {from below and above} by the Bures distance from $\bar{E}$,
\begin{equation}
\setlength{\abovedisplayskip}{3pt}
\setlength{\belowdisplayskip}{3pt}
C_1d_B(\rho,\bar{E}) \leq V(\rho) \leq C_2d_B(\rho,\bar{E}),
\label{C1d<=V<=C2d_2}
\end{equation}
where $C_1 = 1/\sqrt{6}$ and $C_2= 2\sqrt{2}$. It implies
\begin{equation*}
\setlength{\abovedisplayskip}{3pt}
\setlength{\belowdisplayskip}{3pt}
\mathbb{E}(d_B(\rho_t,\bar{E})) \leq \frac{C_2}{C_1}d_B(\rho_0,\bar{E})e^{-\bar{C}t}, \quad \forall \rho_0 \in \mathcal{S},
\end{equation*}
which means that the set $\bar{E}$ is exponentially stable in mean with average Lyapunov exponent less or equal than $-\bar{C}$.

Now we consider the stochastic process
$
Q(\rho_t,t) = e^{\bar{C}t}V(\rho_t) \geq 0
$
whose infinitesimal generator is given by
$
\mathscr{L}Q(\rho,t) = e^{\bar{C}t}( \bar{C}\,V(\rho)+\mathscr{L}V(\rho) )\leq 0.
$ 
Hence, the process $Q(\rho_t,t)$ is a positive supermartingale. Due to Doob's martingale convergence theorem \cite{revuz2013continuous}, the process $Q(\rho_t,t)$ converges almost surely to a finite limit as $t$ tends to infinity. Consequently, $Q(\rho_t,t)$ is almost surely bounded, that is $\sup_{t \geq 0}Q(\rho_t,t) = A$, for some a.s. finite random variable $A$. This implies $\sup_{t \geq 0} V(\rho_t) = Ae^{-\bar{C}t}$ a.s. Letting $t$ goes to infinity, we obtain $\limsup_{t \rightarrow \infty} \frac{1}{t} \log V(\rho_t) \leq -\bar{C}$ a.s. By the inequality~\eqref{C1d<=V<=C2d_2},
\begin{equation} 
\limsup_{t\rightarrow\infty}\frac{1}{t}\log d_B(\rho_t,\bar E) \leq -\bar{C}, \qquad a.s.
\label{rate QSR_2}
\end{equation}
which means that the set $\bar E$ is a.s. exponentially stable with sample Lyapunov exponent less or equal than $-\bar{C}$. 

{Finally, the fact that the probability of convergence to $\bar{\boldsymbol\rho}\in\bar{E}$ is $\mathrm{Tr}(\rho_0 \bar{\boldsymbol\rho})$ may be proved by standard arguments (see e.g~\cite[Theorem~5.1]{liang2019exponential}).}
The proof is complete.\hfill$\square$
\subsection{Exponential stabilization by continuous feedback}
In this {section}, we study the exponential stabilization of system~\eqref{SME} towards a target state $\bar{\boldsymbol\rho}\in\bar{E}$. 
{We first establish} a general result ensuring the exponential convergence towards $\bar{\boldsymbol\rho}$ under some assumptions on the feedback law and an additional local Lyapunov type condition. Next, we design a parametrized family of feedback control laws satisfying such conditions {for some choice of the control Hamiltonian}. {Denote $X_{\bar{\boldsymbol\rho}}(\rho):=\mathrm{Tr}(\rho\bar{\boldsymbol\rho})$ and} $\Theta_u(\rho):=u(\rho)\mathrm{Tr}(i[H_1,\rho]\bar{\boldsymbol\rho})$.
\begin{lemma}
Assume that the initial state satisfies $\rho_0 \neq \bar{\boldsymbol\rho},$ $u\in\mathcal{C}^1(\mathcal{S}\setminus\bar{\boldsymbol\rho},\mathbb{R})$ and 
$|\Theta_u(\rho)|\leq C(1-X_{\bar{\boldsymbol\rho}}(\rho))$
for some $C>0$ . Then
$
\mathbb{P}( \rho_t \neq \bar{\boldsymbol\rho}, \forall\, t\geq 0 )=1.
$
\label{Never reach with feedback 2}
\end{lemma}
\proof
Given $\varepsilon>0,$ we consider any $\mathcal{C}^2$ function on $\mathcal{S}$ such that
\begin{equation*}
\setlength{\abovedisplayskip}{2pt}
\setlength{\belowdisplayskip}{3pt}
V(\rho)=1/(1-X_{\bar{\boldsymbol\rho}}(\rho)),\qquad \mbox{if } X_{\bar{\boldsymbol\rho}}(\rho)<1-\varepsilon.
\end{equation*}
{Under the assumptions of the lemma, it is easy to check that} $\mathscr{L}V(\rho)\leq KV(\rho)$ for some $K>0,$ whenever $X_{\bar{\boldsymbol\rho}}(\rho)<1-\varepsilon$.
To conclude the proof, one just applies the same arguments as in {Lemma~\ref{Never reach without feedback 2_1} and}~\cite[Lemma 4.1]{liang2019exponential}. \hfill$\square$

Generally speaking, {based on the support theorem~\cite{stroock1972support}, trajectories of Equation~\eqref{SME} may be interpreted as limits of solutions of the following deterministic equation}

\begin{equation}
\setlength{\abovedisplayskip}{1pt}
\setlength{\belowdisplayskip}{3pt}
\begin{split}
\dot{\rho}_v(t) = F_0(\rho_v&(t))+\sum^2_{j=1}\widehat{F}_j(\rho_v(t))\sum^2_{j=1}\sqrt{\eta_j}G_j(\rho_v(t))v_j(t),
\end{split}
\label{ODE2}
\end{equation}
with $v_1(t),v_2(t)\in\mathcal{V}$, where $\mathcal{V}$ is the set of all piecewise constant functions from $\mathbb{R}_+$ to $\mathbb{R}$, and 
\begin{equation*}
\setlength{\abovedisplayskip}{3pt}
\setlength{\belowdisplayskip}{3pt}
\begin{split}
\widehat{F}_1(\rho):=&(1-\eta_1)(L_1\rho L_1-\rho)+2\eta_1\mathrm{Tr}(L_1\rho)G_1(\rho),\\
\widehat{F}_2(\rho):=&(1-\eta_2)(L_2\rho L_2-\rho)+2\eta_2\mathrm{Tr}(L_2\rho)G_2(\rho),
\end{split}
\end{equation*}	 
with $F_0$, $G_1$ and $G_2$ defined as in~\eqref{SME}. 
{In particular,} the set $\mathcal S$ is positively invariant for Equation~\eqref{ODE2}. 
\begin{lemma}
Let $\bar{\boldsymbol\rho}=\xi \xi^*$ with $\xi\in\{\Psi_{\pm},\Phi_{\pm}\}.$ Assume that $u\in\mathcal{C}^1(\mathcal{S}\setminus\bar{\boldsymbol\rho},\mathbb{R})$, $u\neq 0$ on the set $\{\rho|\,X_{\bar{\boldsymbol\rho}}(\rho)=0\}$ and $\xi$ is not an eigenvector of $H_1^2$. Then for all $r>0$ and any given initial state $\rho_0 \in \mathcal{S},$
$
\mathbb{P}(\tau_{r} < \infty)=1,
$
where $\tau_{r}: = \inf\{t \geq 0|\, \rho_t \in B_r(\bar{\boldsymbol\rho}) \}$ and $\rho_t$ corresponds to the solution of system~\eqref{SME}.
\label{LEM:ENTER2}
\end{lemma}
\proof
The lemma holds trivially for $\rho_{0} \in B_r(\bar{\boldsymbol\rho})$, as in this case $\tau_{r} = 0$. Let us thus suppose that $\rho_{0} \in \mathcal{S} \setminus B_r(\bar{\boldsymbol\rho})$. We show that there exists $T\in (0,\infty)$ and $\zeta\in (0,1)$ such that $\mathbb{P}_{\rho_0}(\tau_{r}<T )>\zeta$. For this purpose, we make use of the support theorem. Consider the following differential equation derived from~\eqref{ODE2},
\begin{equation*}
\setlength{\abovedisplayskip}{3pt}
\setlength{\belowdisplayskip}{3pt}
\begin{split}
\dot{X}_{\bar{\boldsymbol\rho}}(\rho_{v}(t))=&\Theta_u(\rho_{v}(t))+\sum^2_{j=1}\Theta_j(\rho_{v}(t))+2X_{\bar{\boldsymbol\rho}}(\rho_{v}(t))\sum^2_{j=1}\sqrt{\eta_j M_j}P_j(\rho_{v}(t))v_j(t),
\end{split}
\end{equation*} 
where $v_1(t),v_2(t)\in\mathcal{V}$ is the control input, and
\begin{align*}
\setlength{\abovedisplayskip}{3pt}
\setlength{\belowdisplayskip}{3pt}
\Theta_1(\rho)&:=4\eta_1 M_1\mathrm{Tr}(L_z\rho)P_1(\rho)X_{\bar{\boldsymbol\rho}}(\rho),\quad P_1(\rho):=\bar{\lambda}_z-\mathrm{Tr}(L_z\rho)\text{ where }L_z\bar{\boldsymbol\rho}=\bar{\lambda}_z\bar{\boldsymbol\rho};\\
\Theta_2(\rho)&:=4\eta_2 M_2\mathrm{Tr}(L_x\rho)P_2(\rho)X_{\bar{\boldsymbol\rho}}(\rho),\quad P_2(\rho):=\bar{\lambda}_x-\mathrm{Tr}(L_x\rho)\text{ where }L_x\bar{\boldsymbol\rho}=\bar{\lambda}_x\bar{\boldsymbol\rho}.
\end{align*} 
Denote that $\mathbf{P}_1:=\{\rho\in\mathcal{S}|\,P_1(\rho)=0\}$ and $\mathbf{P}_2:=\{\rho\in\mathcal{S}|\,P_2(\rho)=0\}$, we now analyze the following four different cases ,
\begin{enumerate}
\item If $\bar{\boldsymbol\rho}=\boldsymbol\Psi_{\pm}$, then $\mathbf{P}_1\cap\mathbf{P}_2=\mathcal{S}_{\Psi}\cap\mathcal{R}_{\pm}=\boldsymbol\Psi_{\pm}$;
\item If $\bar{\boldsymbol\rho}=\boldsymbol\Phi_{\pm}$, then $\mathbf{P}_1\cap\mathbf{P}_2=\mathcal{S}_{\Phi}\cap\mathcal{R}_+=\boldsymbol\Phi_{\pm}$.
\end{enumerate}
Suppose $\bar{\boldsymbol\rho}=\boldsymbol\Psi_{+}$. Due to the assumption of the lemma on the feedback law and $H_1$, one can easily show that $X_{\bar{\boldsymbol\rho}}(\rho_v(t))>0$ for $t>0.$ For $t>0,$ we can thus take the feedback $v_1=K P_1(\rho)/X_{\bar{\boldsymbol\rho}}(\rho)$ and $v_2=K P_2(\rho)/X_{\bar{\boldsymbol\rho}}(\rho)$ with $K>0$ sufficiently large.  The proposed control input $v$ guarantees that $\rho_{v}(t)\in B_r(\bar{\boldsymbol\rho})$ for $t\leq T$ with $T<\infty$. The other three Bell states can be treated in the same way. Now, considering the stochastic solution of~\eqref{SME}, we deduce  that $\mathbb P(\rho_t\in B_r(\bar{\boldsymbol\rho}))>0$ for $t\leq T$  from the support theorem~\cite{stroock1972support}.

By compactness of $\mathcal S\setminus B_r(\bar{\boldsymbol\rho})$ and the Feller continuity of $\rho_t,$ we have 
$
\sup_{\rho_0 \in\mathcal S\setminus B_r(\bar{\boldsymbol\rho})} \mathbb{P}_{\rho_0}(\tau_{r}\geq T) \leq 1-\zeta<1,
$ \footnote{Note that $\mathbb{P}_{\rho_0}$ corresponds to the probability law of $\rho_t$ starting at $\rho_0$ and the associated expectation is denoted by $\mathbb{E}_{\rho_0}$.}
for some $\zeta>0.$
By Dynkin inequality~\cite{dynkin1965markov}, 
\begin{equation*}
\setlength{\abovedisplayskip}{3pt}
\setlength{\belowdisplayskip}{3pt}
\sup_{\rho_0 \in\mathcal S\setminus B_r(\bar{\boldsymbol\rho})} \mathbb{E}_{\rho_0}(\tau_{r}) \leq \frac{T}{1-\sup_{\rho_0 \in\mathcal S\setminus B_r(\bar{\boldsymbol\rho}) } \mathbb{P}_{\rho_0}(\tau_{r}\geq T)}\leq \frac{T}{\zeta}.
\end{equation*}
Then by Markov inequality, for all $\rho_0 \in\mathcal S\setminus B_r(\bar{\boldsymbol\rho})$, we have 
\begin{equation*}
\setlength{\abovedisplayskip}{3pt}
\setlength{\belowdisplayskip}{3pt}
\mathbb{P}_{\rho_0}(\tau_{r}=\infty) = \lim_{n\rightarrow \infty} \mathbb{P}_{\rho_0}(\tau_{r} \geq n) \leq \lim_{n\rightarrow \infty} \mathbb{E}_{\rho_0}(\tau_{r})/n=0,
\end{equation*}
which implies
$
\mathbb{P}_{\rho_0}( \tau_{r}<\infty )=1.
$
The proof is  complete.\hfill$\square$

By combining the previous lemmas and following arguments similar to~\cite[Theorem 6.2]{liang2019exponential}, we get the following general result concerning the exponential stabilization towards Bell states. 
\begin{theorem}
Assume that $\rho_0\in \mathcal S$ and the feedback control law satisfies the assumptions of Lemma~\ref{Never reach with feedback 2} and Lemma~\ref{LEM:ENTER2}. Additionally, suppose that there exists a positive-definite function $V(\rho)$ such that $V(\rho)=0$ if and only if $\rho=\bar{\boldsymbol\rho}$, and $V$ is continuous on $\mathcal{S}$ and twice continuously differentiable on the set $\mathcal S\setminus\bar{\boldsymbol\rho}$. Moreover, suppose that there exist positive constants $C$, $C_1$ and $C_2$ such that 
\begin{enumerate}
\item[(i)] $C_1 \, d_B(\rho,\bar{\boldsymbol\rho}) \leq V(\rho) \leq C_2 \, d_B(\rho,\bar{\boldsymbol\rho})$, $\forall\,\rho\in\mathcal{S}$, and 
\item[(ii)] $\limsup_{\rho\rightarrow\bar{\boldsymbol\rho}}\frac{\mathscr{L}V(\rho)}{V(\rho)}=-C$.
\end{enumerate}
Then, $\bar{\boldsymbol\rho}$ is a.s. exponentially stable for the system~\eqref{SME} with sample Lyapunov exponent less or equal than $-C-\frac{K}{2}$, where $K:=\liminf_{\rho \rightarrow\bar{\boldsymbol\rho}}\big(g_1^2(\rho)+g_2^2(\rho)\big)$ with $g_j(\rho):=\sqrt{\eta_j}\frac{\partial V(\rho)}{\partial \rho}\frac{G_j(\rho)}{V(\rho)}$ for $j=1,2$.
\label{Thm a.s. exp stab 2}
\end{theorem}
Next, we derive a general condition on the feedback law and the control Hamiltonian which allows us to apply the previous theorem.
\begin{theorem}
Let $\rho_0\in\mathcal{S}$ and $\bar{\boldsymbol\rho}\in\bar{E}$ be the target state. Suppose that the feedback law and control Hamiltonian satisfy Lemma~\ref{Never reach with feedback 2}, Lemma~\ref{LEM:ENTER2} and the following relation
\begin{equation}
\setlength{\abovedisplayskip}{3pt}
\setlength{\belowdisplayskip}{3pt}
\limsup_{\rho\rightarrow\bar{\boldsymbol\rho}}\Theta_u(\rho)/d_B^2(\rho,\bar{\boldsymbol\rho})=0.
\label{COND_U2}
\end{equation}
Then $\bar{\boldsymbol\rho}$ is almost surely exponentially stable with sample Lyapunov exponent less or equal than $-\bar{C}$ where $\bar{C}=\min\{\eta_1M_1,\eta_2M_2\}.$ 
\label{THM:CONTROL DESING}
\end{theorem}
\proof
To prove the theorem, we show that we can apply Theorem~\ref{Thm a.s. exp stab 2} with the Lyapunov function 
$
V(\rho) = \sqrt{1-X_{\bar{\boldsymbol\rho}}(\rho)}
$
with $\bar{\boldsymbol\rho}\in\bar{E}$. Note that
$
d_B(\rho,\bar{\boldsymbol\rho})\leq V(\rho) \leq \sqrt{2}d_B(\rho,\bar{\boldsymbol\rho}),
$
we are then left to show the condition (ii). The infinitesimal generator of the Lyapunov function satisfies,
%
\begin{equation*}
\setlength{\abovedisplayskip}{3pt}
\setlength{\belowdisplayskip}{3pt}
\mathscr{L}V(\rho)\leq \frac{\Theta_u(\rho)}{2V(\rho)}-\frac{X^2_{\bar{\boldsymbol\rho}}(\rho)\bar{C}}{2V^3(\rho)}\Big((P_1(\rho))^2+(P_2(\rho))^2\Big).
\end{equation*}
Since $\rho\geq 0$, by estimating the right hand side of the above inequality, we obtain the following for all $\rho \in \mathcal S\setminus\bar{\boldsymbol\rho}$, 
\begin{equation*}
\setlength{\abovedisplayskip}{3pt}
\setlength{\belowdisplayskip}{3pt}
\mathscr{L}V(\rho)\leq-\frac{\bar{C}}{2}V(\rho)\left(X^2_{\bar{\boldsymbol\rho}}(\rho)-\frac{\Theta_u(\rho)}{\bar{C}V^2(\rho)}\right).
\end{equation*}
Since $g_1^2(\rho)+g_2^2(\rho)\geq \bar{C}X^2_{\bar{\boldsymbol\rho}}(\rho)$ and by using the relation~\eqref{COND_U2}, we can apply Theorem~\ref{Thm a.s. exp stab 2} with $C=\bar{C}/2$ and $K=\bar C.$ The proof is hence complete.\hfill$\square$

The application of the previous results is given below. 
\begin{proposition}
Consider system~\eqref{SME} with $\rho_0\in\mathcal{S}$. Let $\bar{\boldsymbol\rho}\in\bar{E}$ be the target state. Define the control Hamiltonian as 
$H_1=\sigma_z\otimes\sigma_y-3(\mathds{1}\otimes\sigma_y)$
and the feedback law as 
\begin{equation}
\setlength{\abovedisplayskip}{3pt}
\setlength{\belowdisplayskip}{3pt}
u(\rho)=\alpha (1-X_{\bar{\boldsymbol\rho}}(\rho))^{\beta}-\gamma \mathrm{Tr}(i[H_1,\rho]\bar{\boldsymbol\rho}),
\label{FB:2OM}
\end{equation}
where $\gamma>0,$  $\beta>1$ and $\alpha>0$ sufficiently large. Then $\bar{\boldsymbol\rho}$ is almost surely exponentially stable with sample Lyapunov exponent less or equal than $-\bar{C}$ where $\bar{C}=\min\{\eta_1M_1,\eta_2M_2\}.$ 
\label{THM:DESIGN_2OM}
\end{proposition}
\proof
Since $\rho\geq 0$,
we can show that the feedback law and the control Hamiltonian satisfy the relation~\eqref{COND_U2} and the assumptions of Lemma~\ref{Never reach with feedback 2} and Lemma~\ref{LEM:ENTER2}. The proof is complete by applying Theorem~\ref{THM:CONTROL DESING}. \hfill$\square$
\section{System with one quantum channel}
\label{SEC:1QC}
In this section, our purpose is to  stabilize the system~\eqref{SME} towards a target Bell state $\bar{\boldsymbol\rho}\in\bar{E}$ with only one  measurement operator $L_1=\sqrt{M_1}L_z$. Unlike the previous case, the diffusion terms strengthen the exponential convergence towards $\mathcal{S}_{\Psi}\cup\mathcal{S}_{\Phi}$ instead of $\bar E$. Hence, we do not expect to obtain exponential stabilization towards $\bar{\boldsymbol\rho}$ by using the methods developed before. Hence, we focus on asymptotic stabilization of the system~\eqref{SME}. Here, we take $H_0=\omega L_z$ with $\omega>0$. 


The following result is analogous to Lemma~\ref{LEM:ENTER2}.
\begin{lemma}
Let $\bar{\boldsymbol\rho}=\xi \xi^*$ with $\xi\in\{\Psi_{\pm},\Phi_{\pm}\}.$ Assume that $u\in\mathcal{C}^1(\mathcal{S}\setminus\bar{\boldsymbol\rho},\mathbb{R}^m)$, $u_1\neq 0$ on the set $\{\rho|\,X_{\bar{\boldsymbol\rho}}(\rho)=0\}$ and $\xi$ is not an eigenvector of $H_1^2$. Suppose moreover $u_k\equiv 0$ for $k>1$ on the above set. Then for all $r>0$ and any given initial state $\rho_0 \in \mathcal{S},$
$
\mathbb{P}(\tau_{r} < \infty)=1,
$
where $\tau_{r}: = \inf\{t \geq 0|\, \rho_t \in B_r(\bar{\boldsymbol\rho}) \}$ and $\rho_t$ corresponds to the solution of System~\eqref{SME}.
\label{LEM:ENTER}
\end{lemma}
By employing the first two steps of the proof of~\cite[Theorem 6.2]{liang2019exponential}, we can obtain the general result concerning the asymptotic stabilization of System~\eqref{SME} with only one quantum channel towards the target Bell state.
\begin{theorem}
Assume that the feedback law $u$ satisfies the assumptions of Lemma~\ref{LEM:ENTER}. Additionally, suppose that there exists a positive function $V(\rho)$ such that $V(\rho)=0$ if and only if $\rho=\bar{\boldsymbol\rho}$, and $V$ is continuous on $\mathcal{S}$ and twice continuously differentiable on the set $\mathcal S\setminus\bar{\boldsymbol\rho}$. Moreover, suppose that there exist positive constants $C$, $C_1$ and $C_2$ such that 
\label{THM:AC_1OM}
\begin{enumerate}
\item[(i)] $C_1 \, d^p_B(\rho,\bar{\boldsymbol\rho}) \leq V(\rho) \leq C_2 \, d^p_B(\rho,\bar{\boldsymbol\rho})$ with $p>0$, for all $\rho\in\mathcal{S}$, and 
\item[(ii)] $\mathscr{L}V(\rho)\leq 0,$ for all $\rho \in B_r(\bar{\boldsymbol\rho})$ with $r>0$.
\end{enumerate}
Then, $\bar{\boldsymbol\rho}$ is a.s. asymptotically stable for the system~\eqref{SME}.
\end{theorem}
\begin{remark}
Theorem~\ref{THM:AC_1OM} ensures the global asymptotic stabilization of the system only providing local Lyapunov type condition. The further assumptions on $u_k$ and $H_k$ are used to avoid the presence of invariant subsets of $\mathcal S$. These conditions are not optimal and may be easily weakened. We believe that by applying~\cite[Proposition 4.5]{liang2019exponential}, we can relax these conditions for the case $\eta<1$. We note that  we do not need to find a global Lyapunov condition or apply the LaSalle theorem as in~\cite{mirrahimi2007stabilizing,yamamoto2007feedback}. 
\end{remark}
Next, we define the following continuously differentiable function on $[0,1]$,
%
\begin{equation*}
\setlength{\abovedisplayskip}{1pt}
\setlength{\belowdisplayskip}{1pt}
f(x) = 
\begin{cases}
0,&\text{if }x\in[0,\epsilon);\\
\frac12\sin\big(\frac{\pi(x-1/2)}{1-2\epsilon}\big)+\frac12,&\text{if }x\in[\epsilon,1-\epsilon);\\
1,&\text{if }x\in(1-\epsilon,1],
\end{cases}
\end{equation*}
where $\epsilon\in(0,1/2)$. Denote $\Pi_1:=-\sigma_y\otimes\sigma_z-3(\sigma_y\otimes\mathds{1})$ and $\Pi_2:=\sigma_y\otimes\sigma_z-3(\sigma_y\otimes\mathds{1})$.
Then we propose the following continuous feedback law and control Hamiltonians inspired by~\cite[Theorem 5.1]{mirrahimi2007stabilizing}.
\begin{proposition}
Consider the system~\eqref{SME} with $\rho_0\in\mathcal{S}$. Let $\bar{\boldsymbol\rho}\in\bar{E}$ be the target state. Define $H_1 = \sigma_z\otimes\sigma_y-3(\mathds{1}\otimes\sigma_y)$ and the feedback laws in the following form
\begin{equation}
\setlength{\abovedisplayskip}{3pt}
\setlength{\belowdisplayskip}{3pt}
\begin{split}
u_1(\rho) &= \gamma_1-\mathrm{Tr}(i[H_1,\rho]\bar{\boldsymbol\rho}),\\
u_2(\rho) &= f(X_{\bar{\boldsymbol\rho}}(\rho))(\gamma_2-\mathrm{Tr}(i[H_2,\rho]\bar{\boldsymbol\rho})),
\end{split}
\label{FB:1OM}
\end{equation}
where $|\gamma_1|=|\gamma_2|$ sufficient large. 
If
\begin{itemize}
\item $\bar{\boldsymbol\rho}=\boldsymbol\Psi_{\pm}$, take $\gamma_1=\pm\gamma_2$ and $H_2=\Pi_1$;
\item $\bar{\boldsymbol\rho}=\boldsymbol\Phi_{\pm}$, take $\gamma_1=\pm\gamma_2$ and $H_2=\Pi_2$.
\end{itemize}
Then $\bar{\boldsymbol\rho}$ is a.s. asymptotically stable.
\label{THM:DESIGN_1OM}
\end{proposition}
\proof
We apply Theorem~\ref{THM:AC_1OM} with the Lyapunov function $V(\rho)=1-X_{\bar{\boldsymbol\rho}}(\rho)$. We can easily verify that the feedback law and control Hamiltonians satisfy the assumptions of Lemma~\ref{LEM:ENTER}, $d^2_B(\rho,\bar{\boldsymbol\rho})\leq V(\rho) \leq2d^2_B(\rho,\bar{\boldsymbol\rho})$ in $\mathcal{S}$ and $\mathscr{L}V(\rho)\leq 0$ in a neighborhood of $\bar{\boldsymbol\rho}$. Hence, the proof is complete. \hfill$\square$
\section{Simulations}
\label{SEC:SIM}
In this section, we simulate the dynamics of two-qubit systems in order to illustrate our results. 

\begin{figure}[thpb]
    \includegraphics[width=15cm]{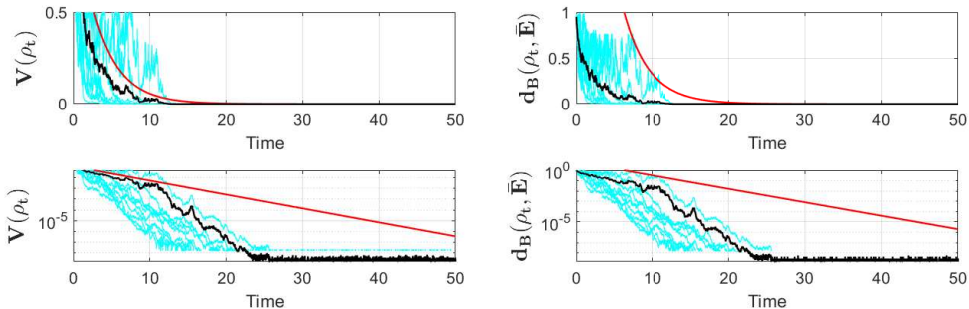}
    \caption{Quantum state reduction with two quantum channels and $u\equiv0$ starting at $diag(0.2,0.3,0.1,0.4)$, when $\omega=0.3$, $\eta_1=0.3$, $M_1=1$, $\eta_2=0.4$, $M_2=0.9$: The black curve represents the mean value of the 10 arbitrary samples, the red curve represents the exponential reference with exponent $-\bar{C}$. The figures at the bottom are the semi-log versions of the ones at the top.}
\label{FIG:QSR_2OM}
\end{figure}

\begin{figure}[thpb]
    \includegraphics[width=15cm]{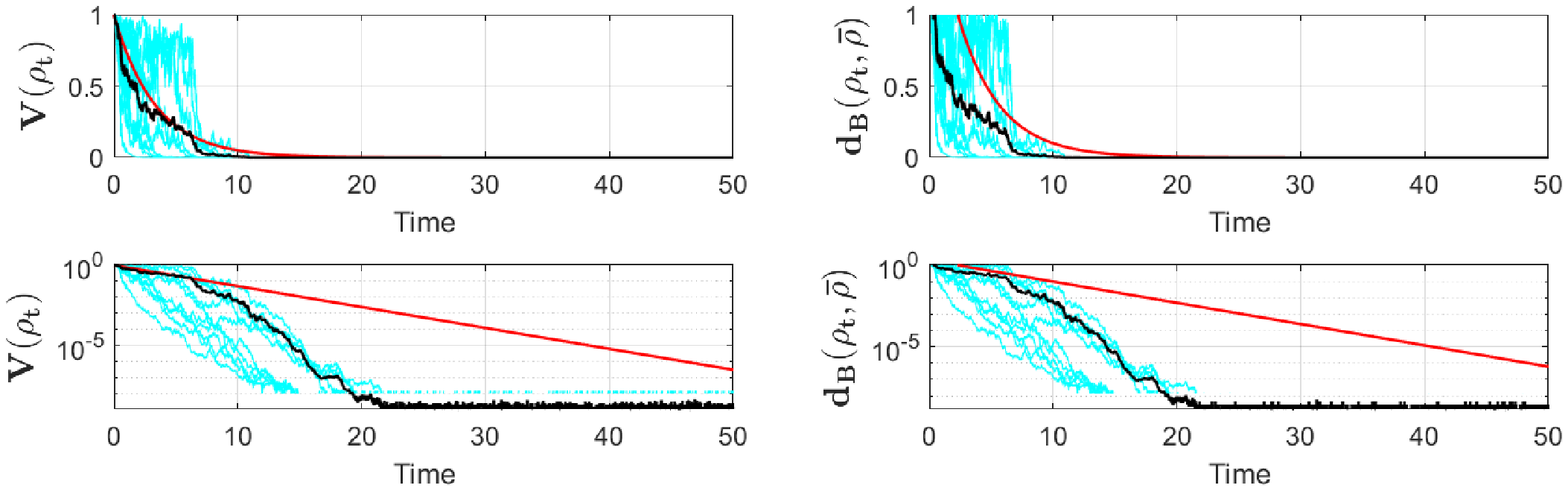}
    \caption{Exponential stabilization of system~\eqref{SME} towards $\boldsymbol\Psi_+$ with the feedback laws~\eqref{FB:2OM} starting at $\boldsymbol\Phi_-$, when $\omega=0.3$, $\eta_1=0.3$, $M_1=1$, $\eta_2=0.4$, $M_2=0.9$, $\alpha=10$, $\beta=12$ and $\gamma=1$: The black curve represents the mean value of the 10 arbitrary samples, the red curve represents the exponential reference with exponent $-\bar{C}$. The figures at the bottom are the semi-log versions of the ones at the top.}
\label{FIG:ES_2OM_SP}
\end{figure}

\begin{figure}[thpb]
    \includegraphics[width=15cm]{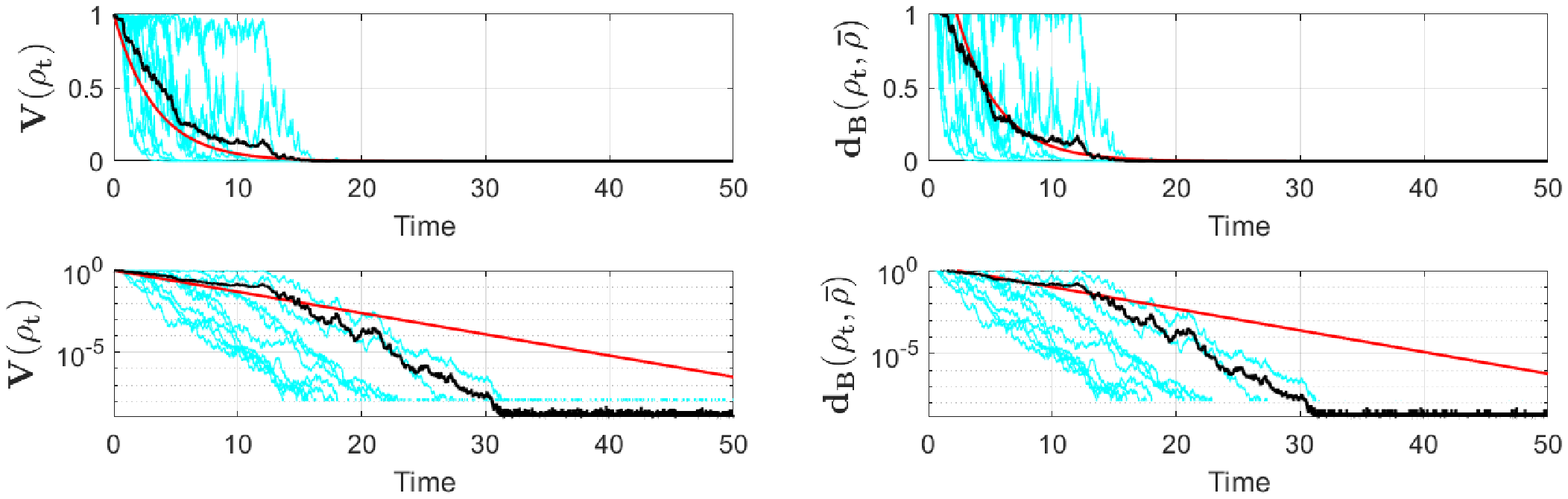}
    \caption{Exponential stabilization of system~\eqref{SME} towards $\boldsymbol\Phi_-$ with the feedback laws~\eqref{FB:2OM} starting at $\boldsymbol\Psi_+$, when $\omega=0.3$, $\eta_1=0.3$, $M_1=1$, $\eta_2=0.4$, $M_2=0.9$, $\alpha=10$, $\beta=12$ and $\gamma=1$: The black curve represents the mean value of the 10 arbitrary samples, the red curve represents the exponential reference with exponent $-\bar{C}$. The figures at the bottom are the semi-log versions of the ones at the top.}
\label{FIG:ES_2OM_HN}
\end{figure}

The simulations in the case with two quantum channels are shown in Fig.~\ref{FIG:QSR_2OM}, Fig.~\ref{FIG:ES_2OM_SP} and
Fig.~\ref{FIG:ES_2OM_HN}. Fig.~\ref{FIG:QSR_2OM} shows the case $u\equiv0$; we observe that the expectation of the Lyapunov function $\mathbb{E}(V(\rho_t))$ is bounded by the exponential function $V(\rho_0)e^{-\bar{C}t}$, and the expectation of the Bures distance $\mathbb{E}(d_B(\rho_t,\bar{E}))$ is always bounded by $4\sqrt{3}d_B(\rho_0,\bar{E})e^{-\bar{C}t}$ which confirms the results of Theorem~\ref{QSR_2}. Then we set $\boldsymbol\Psi_+$ as the target state and $\boldsymbol\Phi_-$ as the initial state; the behavior of the system with the continuous feedback~\eqref{FB:2OM} is shown in Fig.~\ref{FIG:ES_2OM_SP}. 
Similar simulations for the case with $\boldsymbol\Phi_-$ as the target state and $\boldsymbol\Psi_+$ as initial state are shown in Fig.~\ref{FIG:ES_2OM_HN}.

\begin{figure}[thpb]
    \includegraphics[width=15cm]{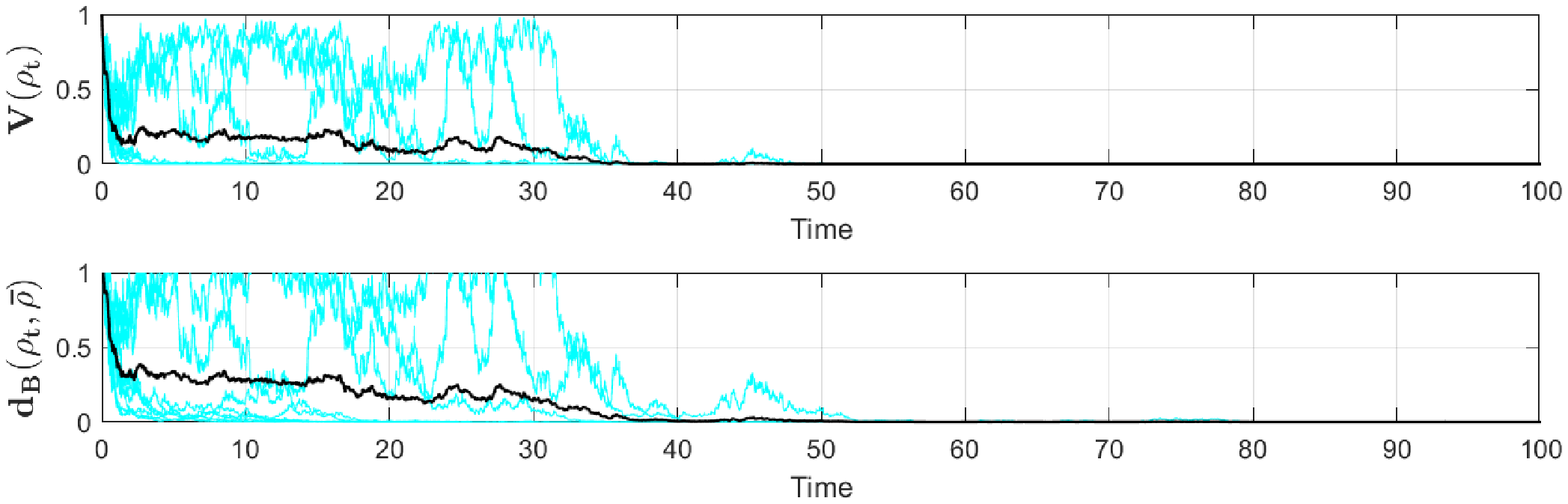}
    \caption{Asymptotic stabilization of system~\eqref{SME} towards $\boldsymbol\Psi_+$ with the feedback laws~\eqref{FB:1OM} starting at $\boldsymbol\Phi_-$, when $\omega=0.3$, $\eta_1=0.3$, $M_1=1$, $\epsilon=0.15$ and $\gamma_1=\gamma_2=4$: The black curve represents the mean value of the 10 arbitrary samples.}
\label{FIG:AS_1OM_SP}
\end{figure}

\begin{figure}[thpb]
    \includegraphics[width=15cm]{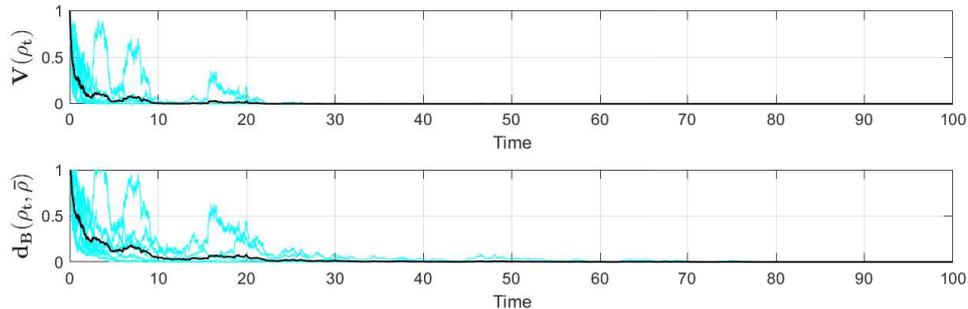}
    \caption{Asymptotic stabilization of system~\eqref{SME} towards $\boldsymbol\Phi_-$ with the feedback laws~\eqref{FB:1OM} starting at $\boldsymbol\Psi_+$, when $\omega=0.3$, $\eta_1=0.3$, $M_1=1$, $\epsilon=0.15$ and $\gamma_1=\gamma_2=4$: The black curve represents the mean value of the 10 arbitrary samples.}
\label{FIG:AS_1OM_HN}
\end{figure}

The simulations in the case with only one quantum channel are shown in Fig.~\ref{FIG:AS_1OM_SP} and
Fig.~\ref{FIG:AS_1OM_HN} for $\boldsymbol\Psi_+$ as the target state and $\boldsymbol\Phi_-$ as the target state respectively. Such simulations clearly confirm the validity of Proposition~\ref{THM:DESIGN_1OM}.
\section{Conclusion and future works}
In this paper, we have studied the asymptotic behavior of trajectories of two-qubit systems. In particular, for the case of two quantum channels, we have provided a general result concerning the feedback exponential stabilization towards the Bell states by applying local stochastic Lyapunov techniques and analyzing the asymptotic behavior of quantum trajectories. Furthermore, we constructed a parameterized continuous feedback law satisfying the conditions of our general results. Next, for the system with one quantum channel, by a similar analysis, we proposed a continuous feedback law stabilizing the system asymptotically. 

Further research lines will address the possibility of extending our results to multi-qubit entanglement generation, or in presence of delays.



\bibliographystyle{unsrt}
\bibliography{LIANG_CDC19}

\begin{thebibliography}{10}

\bibitem{nielsen2010quantum}
M.~A. Nielsen and I.~L. Chuang.
\newblock {\em Quantum computation and quantum information}.
\newblock Cambridge university press, 2010.

\bibitem{bengtsson2017geometry}
I.~Bengtsson and K.~{\.Z}yczkowski.
\newblock {\em Geometry of quantum states: an introduction to quantum
  entanglement}.
\newblock Cambridge University Press, 2017.

\bibitem{van2005feedback}
R.~Van~Handel, J.~K. Stockton, and H.~Mabuchi.
\newblock Feedback control of quantum state reduction.
\newblock {\em IEEE Transactions on Automatic Control}, 50(6):768--780, 2005.

\bibitem{armen2002adaptive}
M.~A. Armen, J.~K. Au, J.~K. Stockton, A.~C. Doherty, and H.~Mabuchi.
\newblock Adaptive homodyne measurement of optical phase.
\newblock {\em Physical Review Letters}, 89(13):133602, 2002.

\bibitem{mirrahimi2007stabilizing}
M.~Mirrahimi and R.~Van~Handel.
\newblock Stabilizing feedback controls for quantum systems.
\newblock {\em SIAM Journal on Control and Optimization}, 46(2):445--467, 2007.

\bibitem{tsumura2008global}
K.~Tsumura.
\newblock Global stabilization at arbitrary eigenstates of n-dimensional
  quantum spin systems via continuous feedback.
\newblock In {\em American Control Conference, 2008}, pages 4148--4153, 2008.

\bibitem{ahn2002continuous}
C.~Ahn, A.~C. Doherty, and A.~J. Landahl.
\newblock Continuous quantum error correction via quantum feedback control.
\newblock {\em Physical Review A}, 65(4):042301, 2002.

\bibitem{yamamoto2007feedback}
N.~Yamamoto, K.~Tsumura, and S.~Hara.
\newblock Feedback control of quantum entanglement in a two-spin system.
\newblock {\em Automatica}, 43(6):981--992, 2007.

\bibitem{mabuchi2005principles}
H.~Mabuchi and N.~Khaneja.
\newblock Principles and applications of control in quantum systems.
\newblock {\em International Journal of Robust and Nonlinear Control:
  IFAC-Affiliated Journal}, 15(15):647--667, 2005.

\bibitem{belavkin1983theory}
V.~P. Belavkin.
\newblock On the theory of controlling observable quantum systems.
\newblock {\em Avtomatika i Telemekhanika}, (2):50--63, 1983.

\bibitem{liang2018exponential}
W.~Liang, N.~H. Amini, and P~Mason.
\newblock On exponential stabilization of spin-$\frac12$ systems.
\newblock In {\em IEEE Conference on Decision and Control}, pages 6602--6607,
  2018.

\bibitem{liang2019exponential}
W.~Liang, N.~H. Amini, and P.~Mason.
\newblock On exponential stabilization of {$N$}-level quantum angular momentum
  systems.
\newblock {\em arXiv preprint arXiv:1902.05879}, 2019.
\newblock Submitted.

\bibitem{hudson1984quantum}
R.~L. Hudson and K.~R. Parthasarathy.
\newblock Quantum {I}to's formula and stochastic evolutions.
\newblock {\em Communications in Mathematical Physics}, 93(3):301--323, 1984.

\bibitem{belavkin1989nondemolition}
V.~P. Belavkin.
\newblock Nondemolition measurements, nonlinear filtering and dynamic
  programming of quantum stochastic processes.
\newblock In {\em Modeling and Control of Systems}, pages 245--265. Springer,
  1989.

\bibitem{bouten2007introduction}
L.~Bouten, R.~Van~Handel, and M.~R. James.
\newblock An introduction to quantum filtering.
\newblock {\em SIAM Journal on Control and Optimization}, 46(6):2199--2241,
  2007.

\bibitem{mao2007stochastic}
X.~Mao.
\newblock {\em Stochastic differential equations and applications}.
\newblock Woodhead Publishing, 2007.

\bibitem{khasminskii2011stochastic}
R.~Khasminskii.
\newblock {\em Stochastic stability of differential equations}, volume~66.
\newblock Springer, 2011.

\bibitem{abe2008analysis}
T.~Abe, T.~Sasaki, S.~Hara, and K.~Tsumura.
\newblock Analysis on behaviors of controlled quantum systems via quantum
  entropy.
\newblock {\em IFAC Proceedings Volumes}, 41(2):3695--3700, 2008.

\bibitem{vu2012real}
T.~L. Vu, S.~S. Ge, and C.~C. Hang.
\newblock Real-time deterministic generation of maximally entangled two-qubit
  and three-qubit states via bang-bang control.
\newblock {\em Physical Review A}, 85(1):012332, 2012.

\bibitem{revuz2013continuous}
D.~Revuz and M.~Yor.
\newblock {\em Continuous martingales and Brownian motion}, volume 293.
\newblock Springer, 2013.

\bibitem{stroock1972support}
D.~W. Stroock and S.~R. Varadhan.
\newblock On the support of diffusion processes with applications to the strong
  maximum principle.
\newblock In {\em Proceedings of the Sixth Berkeley Symposium on Mathematical
  Statistics and Probability (Univ. California, Berkeley, Calif., 1970/1971)},
  volume~3, pages 333--359, 1972.

\bibitem{dynkin1965markov}
E.~B. Dynkin.
\newblock Markov processes.
\newblock In {\em Markov Processes}, pages 77--104. Springer, 1965.

\end{thebibliography}
\end{document}